\documentclass[11pt]{article}

\usepackage{amscd,amsmath, amssymb, fancyhdr, epsfig, color, tocloft, mathtools}
\usepackage{dutchcal}

\usepackage[backref=page]{hyperref}
\renewcommand*{\backrefalt}[4]{%
	\ifcase #1 (Not cited.)%
	\or        (Cited on page~#2.)%
	\else      (Cited on pages~#2.)%
	\fi}

\hypersetup{
	colorlinks   = true,
	citecolor    = magenta,
	linkcolor    = blue,
	urlcolor     = magenta	
}

\newcommand{\version}{version 1.0,\ \ February 22, 2022}

\makeatletter
\def\x@arrow{\DOTSB\Relbar}
\def\xlongequalsignfill@{\arrowfill@\x@arrow\Relbar\x@arrow}
\providecommand{\xlongequal}[2][]{%
	\ext@arrow 0099\xlongequalsignfill@{#1}{#2}}
\def\xlongrightarrowfill@{\arrowfill@\relbar\relbar\longrightarrow}
\newcommand{\xlongrightarrow}[2][]{%
	\ext@arrow 0099\xlongrightarrowfill@{#1}{#2}}
\makeatother

\numberwithin{equation}{section}

\def\eqref#1{(\ref{#1})}

\newcommand{\Z}{{\mathbb Z}}
\newcommand{\C}{{\mathbb C}}

\def\1{\sqrt{-1}\,}

\newcommand{\cntrct}                % contraction with a vector field
{\hspace{2pt}\raisebox{1pt}{\text{$\lrcorner$}}\hspace{2pt}}

\newcommand{\arrow}{{\:\longrightarrow\:}}

\newcommand{\calo}{{\cal O}}

\renewcommand{\bar}{\overline}
\renewcommand{\phi}{\varphi}
\renewcommand{\epsilon}{\varepsilon}
\renewcommand{\geq}{\geqslant}
\renewcommand{\leq}{\leqslant}

\newcommand{\im}{\operatorname{im}}
\newcommand{\End}{\operatorname{End}}

\newcommand{\Id}{\operatorname{Id}}

\renewcommand{\sup}{{\operatorname{{\sf sup}}}}

\renewcommand{\dim}{\operatorname{\sf dim}}

\newcommand{\Spec}{\operatorname{{\sf Spec}}}

\newcounter{Mycounter}[section]
\newcounter{lemma}[section]
\newcounter{claim}[section]
\newcounter{sublemma}[section]
\newcounter{corollary}[section]
\renewcommand{\thecorollary}{{Corollary \thesection.\arabic{corollary}}}
\newcommand{\corollary}{%
	\setcounter{corollary}{\value{Mycounter}}
	\refstepcounter{corollary}
	\stepcounter{Mycounter}
	{\noindent \bf \thecorollary:\ }}

\newcounter{theorem}[section]
\renewcommand{\thetheorem}{{Theorem \thesection.\arabic{theorem}}}
\newcommand{\theorem}{%
	\setcounter{theorem}{\value{Mycounter}}
	\refstepcounter{theorem}
	\stepcounter{Mycounter}
	{\noindent \bf \thetheorem:\ }}

\newcounter{conjecture}[section]
\newcounter{proposition}[section]
\newcounter{definition}[section]
\renewcommand{\thedefinition} {{Definition~\thesection.\arabic{definition}}}
\newcommand{\definition}{%
	\setcounter{definition}{\value{Mycounter}}
	\refstepcounter{definition}
	\stepcounter{Mycounter}
	{\noindent \bf \thedefinition:\ }}

\newcounter{example}[section]
\renewcommand{\theexample}{{Example \thesection.\arabic{example}}}
\newcommand{\example}{%
	\setcounter{example}{\value{Mycounter}}
	\refstepcounter{example}
	\stepcounter{Mycounter}
	{\noindent \bf \theexample:\ }}

\newcounter{remark}[section]
\renewcommand{\theremark}{{Remark \thesection.\arabic{remark}}}
\newcommand{\remark}{%
	\setcounter{remark}{\value{Mycounter}}
	\refstepcounter{remark}
	\stepcounter{Mycounter}
	{\noindent \bf \theremark:\ }}

\newcounter{problem}[section]
\newcounter{question}[section]
\makeatletter

\def\blacksquare{\hbox{\vrule width 5pt height 5pt depth 0pt}}
\def\endproof{\blacksquare}

\newcommand{\proof}{{\bf Proof: \ }}
\newcommand{\pstep}{{\bf Proof. Step 1: \ }}

\begin{document}
		
		\begin{center}
			{\Large\bf  Non-linear Hopf
                          manifolds \\[3mm] are locally conformally K\"ahler}\\[5mm]
			%%%%%%%%%%%%%%%%%%%%%%%%%%%%%%%%%%%%%%%%%%%%%%%%%%%%%%%%%%%%
			{\large
				Liviu Ornea\footnote{Liviu Ornea is  partially supported by Romanian Ministry of Education and Research, Program PN-III, Project number PN-III-P4-ID-PCE-2020-0025, Contract  30/04.02.2021},  
				Misha Verbitsky\footnote{Misha Verbitsky is partially supported by
					by the HSE University Basic Research Program, FAPERJ E-26/202.912/2018 
					and CNPq - Process 313608/2017-2.\\[1mm]
					
				}\\[4mm]
				
			}
			
		\end{center}

		{\small
			\hspace{0.15\linewidth}
			\begin{minipage}[t]{0.7\linewidth}
				{\bf Abstract} \\ 
A  Hopf manifold is a quotient of $\C^n\backslash
0$ by the cyclic group generated by a holomorphic
contraction. Hopf manifolds are diffeomorphic to
$S^1\times S^{2n-1}$ and hence do not admit K\"ahler
metrics. It is known that Hopf manifolds defined by linear
contractions (called linear Hopf manifolds)  have locally
conformally K\"ahler (LCK) metrics. In this paper we prove
that the Hopf manifolds defined by non-linear holomorphic
contractions admit 
holomorphic embeddings into linear Hopf manifolds, and,
moreover they admit LCK metrics.
			\end{minipage}
		}
		%%%%%%%%%%%%%%%%%%%%%%%%%%%%%%%%%%%%%%%%%%%%%%%%

	\tableofcontents

%%%%%%%%%%%%%%%%%%%%%%%%%%%%%%%%%%%%%%%%%%%%%%%%%%%%%%%%%%%%	
\section{Introduction}
%%%%%%%%%%%%%%%%%%%%%%%%%%%%%%%%%%%%%%%%%%%%%%%%%%%%%%%%%%%%

A {\bf Hopf manifold} is a quotient of $\C^n \backslash 0$ by a cyclic group $\langle A \rangle$ 
generated by a holomorphic contraction $A$. This definition is motivated by Kodaira's
research on complex surfaces. In \cite[Theorem 34]{_Kodaira_Structure_II_}
Kodaira proved that a complex surface with $b_1=1$ and
$b_2=0$ is a Hopf surface, if it 
contains a curve and has no non-constant meromorphic
functions. Using the Poincar\'e classification
of contractions of $\C^n$, and the further works by 
Dulac, Latt\`es and Sternberg 
\cite{_Lattes_,_Dulac_,_Sternberg_contraction_},
Kodaira has written a normal form 
for a Hopf surface, expressing a given
contraction as a polynomial map
in appropriate complex coordinates.

If the contraction $A$ is linear, then we speak about a {\bf linear Hopf manifold}.

It is easily seen that a Hopf manifold (be it linear or
non-linear) of complex dimension $n$ is  diffeomorphic
with $S^1\times S^{2n-1}$ and hence does not admit
K\"ahler metrics. It is then natural to ask what kind of
Hermitian metric can a Hopf manifold bear. For linear Hopf
surfaces, the answer is known since long: they have
locally conformally K\"ahler (LCK) metrics (see
Section \ref{_LCK_background_subsection_}, also \cite{do}). This
was proven by various authors, first for diagonal linear
contractions (\cite{va_isr}), then for non-diagonalizable
linear contractions (\cite{go, bel}). In \cite{kor},
this constructon was generalized to arbitrary dimension.
This construction of LCK metric was based on deformation
arguments. First, one builds an LCK metric on a linear Hopf
manifold given by diagonalizable linear contrations
(we call such manifolds ``diagonal Hopf manifolds''),
and then one extends this metric to the non-diagonal
linear Hopf manifolds using deformation stability of
LCK metric with potential (\ref{_LCK_pote_def_},
\cite{ov_lckpot}).

In \cite{ov_pams}, we gave a more conceptual construction of
LCK potential on any linear Hopf manifold.

%Moreover, we proved
%that the property of admitting such metrics is hereditary:
%all closed complex submanifolds of an LCK manifold with
%potential are LCK with potential
%(\ref{_hereditary_LCK_pot_remark_}). 

%Following Poincar\'e, Dulac and others, 
%Kodaira made a distinction between 
%{\bf resonant} and {\bf non-resonant}\index{manifold!Hopf!(non-)resonant}
%Hopf manifolds. Let $H= \frac{\C^n \backslash 0}{\langle \gamma \rangle}$
%be a Hopf manifold. Using Hartogs theorem, we extend\index{theorem!Hartogs}
%$\gamma$ to a holomorphic automorphism of $\C^n$; let
%$\alpha_1, ..., \alpha_n$ be the eigenvalues of 
%the differential of $\gamma$ in 0.
%
%A Hopf manifold is
%called {\bf non-resonant} if the numbers
%$\alpha_i$ don't satisfy a relation
%of form $\alpha_i = \prod_{j=1}^d {\alpha_{l_j}}$
%for any $d>1$. Kodaira has shown that
%all non-resonant Hopf surfaces are linear,
%and classified the resonant ones using the
%normal form provided by the Sternberg's 
%version of the Latt\`es-Poincar\'e-Dulac
%theorem.
%
%We re-prove a part of this classification theorem in 
%Section \ref{_flat_Hopf_Section_}, showing that 
%any non-resonant Hopf manifold is linear.
%The full strength of Sternberg (and Poincar\'e-Dulac)
%classification can be found in \cite{_Arnold:ODE+_}
%and \cite{_Sternberg_contraction_}; its application
%to the Hopf manifold is straightforward.

In this paper, we complete this picture by proving that also non-linear Hopf manifolds
do bear  LCK metrics with potential (\ref{_hopf_are_lck_corollary_}).

One way to attack the problem is to find a normal form for a general holomorphic contraction. In dimension greater that 2, this requires a careful analysis of resonant contractions (a case not considered in \cite{_Sternberg_contraction_}). This is the method suggested by Belgun in \cite{_belgun_cirm_}.

We didn't follow this path. Our proof consists in 
 embedding a non-linear Hopf manifold into a  linear Hopf manifold. 

In fact, the result we prove is even stronger.
Let $V$ be a Stein variety with an isolated
singularity $x$, and $\gamma:\; V \arrow V$ an
invertible holomorphic contraction centered in $x\in V$.
Then $M:=\frac{V\backslash x}{\langle \gamma\rangle}$
is a compact, smooth complex variety. We prove that $M$ can be
embedded to a linear Hopf manifold. Then 
\ref{_hereditary_LCK_pot_remark_} implies that
$M$ is LCK with potential whenever it is smooth.

The proof of this result relies on several results in
functional analysis, mainly on the Riesz-Schauder theorem,
that we briefly recall in Subsections
\ref{_RS_subsection_} and
\ref{_normal_family_subsection_}. We used a similar
approach in  \cite{ov_indam}.

%%%%%%%%%%%%%%%%%%%%%%%%%%%%%%%%%%%%%%%%%%%%%%%%%%%%%%%%%%
\section{Prerequisites}
%%%%%%%%%%%%%%%%%%%%%%%%%%%%%%%%%%%%%%%%%%%%%%%%%%%%%%%%%%

%%%%%%%%%%%%%%%%%%%%%%%%%%%%%%%%%%%%%%%%%%%%%%%%%%%%%%%%%%%%
\subsection{Locally conformally K\"ahler manifolds}\label{_LCK_background_subsection_}
%%%%%%%%%%%%%%%%%%%%%%%%%%%%%%%%%%%%%%%%%%%%%%%%%%%%%%%%%%%%%

We gather here the minimal necessary background on locally
conformally K\"ahler structures, as is needed for our main
result. 

\hfill

\definition
A complex Hermitian manifold $(M, I)$ is {\bf locally conformally
K\"ahler} (LCK, for short) if it admits a
covering $(\tilde M, I)$ equipped with a K\"ahler metric
$\tilde\omega$ such that the deck group of the
cover acts on $(\tilde M, \tilde \omega)$ by holomorphic 
homotheties. An {\bf LCK metric} on an LCK manifold
is an Hermitian metric on $(M,I)$ such that its
pullback to $\tilde M$ is conformal with $\tilde \omega$.

\hfill

A special class of LCK manifolds is the following:

\hfill

\definition\label{_LCK_pote_def_} An LCK manifold has  {\bf LCK potential} if it
admits a K\"ahler covering on which the K\"ahler metric
has a global and positive  potential function $\psi$ such that
the deck group multiplies $\psi$ by a
constant.\footnote{In the sequel, a
  differential form which is multiplied by a constant
  factor by the action of the deck group is called {\bf
    automorphic}.} 
In this case, $M$ is called {\bf LCK manifold with
	potential}.

\hfill

\example Let  $H:=(\C^n\setminus 0)/\langle A\rangle$, $A=\lambda\Id$, 
 $|\lambda|> 1$, be a {\bf diagonal Hopf manifold}. It is LCK because the flat K\"ahler metric $\tilde g_0=\sum dz_i\otimes d\bar z_i$ is multiplied by 2 by the deck group $\Z$; and it is LCK with potential because $\tilde g_0$ has the global automorphic potential $\psi:=\sum |z_i|^2$. 
 
 By \cite{ov_lckpot}, compact LCK manifolds with potential
 are stable to small deformations of the complex
 structure. This implies that
 all {\bf linear Hopf manifolds} $(\C^n\setminus 0)/\langle A\rangle$, $A\in\mathrm{GL}(n,\C)$
 with eigenvalues of absolute value $> 1$, are LCK with potential.

\hfill

\remark By \cite{ovv1}, the blowup of an LCK manifold
along a K\"ahler submanifold is still LCK. However, its
K\"ahler cover cannot admit global potentials. Hence, the
class of LCK manifolds with potential is strict.
% Indeed, the LCK Inoue surfaces (\cite{tric, bel})and the LCK Oeljeklaus-Toma manifolds (\cite{ot}) cannot have any LCK metric with potential, see e.g. \cite{oti2}. 

\hfill

\remark\label{_hereditary_LCK_pot_remark_} Clearly, any
complex submanifold of an LCK manifold with potential is
also LCK with potential (\cite{ov_lckpot}).

\hfill

\remark Let $M$ be a compact LCK manifold
with potential. We say that the potential is {\bf proper}
if it is proper as a smooth function, that is, all its
level sets are compact.
It was shown in \cite{ov_jgp_09, ov_jgp_16} that  
 the LCK metric can be approximated 
 in the ${C}^\infty$-topology by an LCK
metric with {\em proper} LCK potential on the same
complex manifold. Therefore, as long as we are interested
in the complex geometry of LCK manifolds with potential,
we may assume that the potential  is a proper
function. This, in turn, is equivalent with the deck group
of the K\"ahler covering being isomorphic to $\Z$
(\cite{ov_lckpot}).

\hfill

\theorem {(\cite{ov_lckpot,ov_pams})}\label{potcon}
Let $M$ be a compact LCK manifold with proper potential, $\dim_\C M\geq 3$, 
and $\tilde M$ its K\"ahler $\Z$-covering.
Then the metric completion $\tilde M_c$
admits a structure of a complex variety, 
compatible with the complex structure on
$\tilde M \subset\tilde M_c$, and 
the complement $\tilde M_c\setminus  \tilde M$
is just one point. Moreover, $\tilde M_c$ is an affine
algebraic variety obtained as an affine cone
over a projective orbifold.\footnote{$\tilde M_c$ is indeed the {\bf Stein completion} of $\tilde M$ in the sense of \cite{andreotti_siu}. The restriction $\dim_\C M\geq 3$ is imposed by our proof which makes use of the filling theorem 
	by Rossi and Andreotti-Siu (\cite{rossi,andreotti_siu}).} 

\hfill

This result is essential in
the proof of the following Kodaira 
embedding type result:

\hfill

\theorem\label{_Embedding_LCK_pot_in_Hopf_}(\cite{ov_lckpot})
Any compact LCK manifold with potential, of complex dimension at least 3, admits a
holomorphic embedding into a linear Hopf manifold.

%%%%%%%%%%%%%%%%%%%%%%%%%%%%%%%%%%%%%%%%%%%%%%%%%%%%%%%
\subsection{The Riesz-Schauder theorem}\label{_RS_subsection_}
%%%%%%%%%%%%%%%%%%%%%%%%%%%%%%%%%%%%%%%%%%%%%%%%%%%%%

The following theorem can be used to obtain a version of Jordan normal
form for a compact operator on a Banach space. Recall
that the {\bf spectrum} of a linear operator $F$ is the set of all
$\mu\in \C$ such that $F - \mu\Id$ is not invertible.

\hfill

\theorem{ (Riesz-Schauder, \cite[Section 5.2]{friedman})}\\
\label{_Riesz_Schauder_main_Theorem_}
Let $F:\; V \arrow V$ be a compact operator on a Banach
space.  Then the spectrum $\Spec F\subset \C$ is 
compact and discrete outside of $0\in \C$. Moreover, 
 for each non-zero $\mu \in \Spec F$, there exists a sufficiently
big number $N\in \Z^{\geq 0}$ such that for each $n>N$ {one has 
\[ V= \ker(F-\mu\Id)^n \oplus \overline{\im (F-\mu\Id)^n},
\]
	where $\overline{\im (F-\mu\Id)^n}$ is the closure of the image.}
Finally, the space $\ker(F-\mu\Id)^n$ is finite-dimensional.

\hfill

\remark\label{_root_space_RS_Remark_}
Recall that   {\bf the root space of an operator \index{root space}
	$F\in \End(V)$, associated with an eigenvalue $\mu$}, 
is $\bigcup_{n\in \Z} \ker(F-\mu\Id)^n$. 
In the finite-dimensional case, the root spaces
coincide with the  Jordan cells of the corresponding 
matrix. Then 
\ref{_Riesz_Schauder_main_Theorem_} can be reformulated
by saying that any compact operator $F\in \End(V)$ admits a 
Jordan cell decomposition, with $V$ identified with\index{Jordan cell}
a completed direct sum of the root spaces, which are
all finite-dimensional; moreover, the eigenvalues
of $F$ converge to zero.

\hfill

We shall need the following corollary of the Riesz-Schauder theorem, which is
obtained using the same arguments as in the finite-dimensional case.

\hfill

\theorem\label{rs}
Let  $F:\; V \arrow V$ be a compact operator on a Banach space.
We say that $v$ is {\bf a root vector}\index{vector!root} for $F$ if $v$ lies in 
a root space of $F$, for some eigenvalue $\mu\in \C$. 
Then the space generated by the root vectors is dense in $V$.
\endproof

\subsection{Normal families and $F$-finite vectors}\label{_normal_family_subsection_}

\definition \label{normal_family}
Let $M$ be a complex manifold, 
and ${\cal F}\subset H^0(\calo_M)$ a 
family of holomorphic functions. We call ${\cal F}$ 
{\bf a normal family} if for each compact
$K\subset M$ there exists $C_K>0$ such that
for each $f\in {\cal F}$, $\sup_K |f| \leq C_K$.

\hfill

\definition 
The {\bf $C^0$-topology} on the space
of functions on $M$ is the topology
of uniform convergence on compacts.

\hfill

\theorem \label{_Montel_Theorem_}
(Montel theorem)
Let ${\cal F}\subset H^0(\calo_M)$
be a normal family of functions,
and $\bar {\cal F}\subset H^0(\calo_M)$ its closure
in the $C^0$-topology. Then $\bar {\cal F}$ is compact
in $C^0$-topology.\footnote{A set which has compact
  closure is called {\bf precompact}.}

\proof 
\cite{_Wu:Montel_}. \endproof

\hfill

\definition\label{_F_finite_Definition_}
Let $F\in \End(V)$ be an endomorphism of a vector space.
A vector $v\in V$ is called {\bf $F$-finite}
if the space generated by $v, F(v), F(F(v)),  ...$
is finite-dimensional. 

\hfill

\remark Clearly, a vector is 
$F$-finite if and only if it is obtained as a combination
of root vectors. Then  \ref{rs} 
implies the following.

\hfill

\corollary\label{_finite_RS_Corollary_}
Let  $F:\; V \arrow V$ be a compact operator on a Banach space,
and $V_0\subset V$ the space of all $F$-finite vectors.
Then $V_0$ is dense in $V$.
\endproof

\hfill

%\remark  \label{_Riesz_theorem_Remark_}
%By Riesz theorem (\cite[Chapter 1, Theorem 4]{_Diestel_}), a closed ball in a normed
%vector space $V$ is never compact, unless $V$ is
%finite-di\-men\-si\-o\-nal. This means that $(H^0(\calo_M),  C^0)$
%does not admit a norm.  A topological vector space in which any bounded subset is
%precompact is called a {\bf Montel space}.

%%%%%%%%%%%%%%%%%%%%%%%%%%%%%%%%%%%%%%%%%%%%%%%%%%%%%%%%%%%%
\section{Holomorphic contractions on Stein varieties}
%%%%%%%%%%%%%%%%%%%%%%%%%%%%%%%%%%%%%%%%%%%%%%%%%%%%%%%%%%%%

Recall that {\bf a holomorphic contraction}
of a variety $V$ with center in $x\in V$
is a map $\gamma:\; V \arrow V$, $\gamma(x)=x$, such that
for each compact subset $K \subset V$ and
any open neighbourhood $U \ni x$, there
exists $N >0$ such that $\gamma^n(K) \subset U$
for any integer $n>N$. In this case, 
$\gamma$ is called {\bf the contraction centered in $x$}.

\hfill

%\ref{_Stein_by_contract_to_linear_Hopf_Theorem_} below
%is a generalization, in a sense, of \ref{embedding},
%and their proofs are very similar. 
%We decided to include both theorems
%because the argument proving \ref{_Stein_by_contract_to_linear_Hopf_Theorem_}
%is more abstract and harder to digest. 
%Moreover, the embedding for LCK manifolds
%with potential has a cohomological interpretation
%(\ref{_embedding_to_Hopf_via_BC_class_Corollary_})
%which the seemingly more general 
%\ref{_Stein_by_contract_to_linear_Hopf_Theorem_}
%lacks. Therefore it
%makes sense to discuss these two theorems separately.

\hfill

%%%%%%%%%%%%%%%%%%%%%%%%%%%%%%%%%%%%%%%%%%%%%%%%
\theorem\label{_Stein_by_contract_to_linear_Hopf_Theorem_}
Let $V$ be a Stein variety, and $\gamma:\; V \arrow V$ 
an invertible holomorphic contraction centered in $x$.
Assume that the group $\langle \gamma\rangle$ acts on $V\backslash x$ properly
discontinuously. Then the quotient variety 
$M:=\frac{V\backslash x}{\langle \gamma\rangle}$ is
compact. Moreover, it admits
a holomorphic embedding to a linear Hopf manifold. 

\hfill

\pstep
Let $V_0 \ni x$ be a precompact Stein neighbourhood of $x\in V$.
By the definition of a contraction, for $k$ sufficiently big,
we have $\gamma^k(V_0) \Subset V_0$.\footnote{Here, as elsewhere,
	the notation $A \Subset B$ means that $A$ is {\em relatively compact} in $B$,
	that is, $A$ is a subset of $B$ and its closure is compact.}
Let $H^0_b(V_0, \calo_V)$ be the space of bounded holomorphic
functions on $V_0$, equipped with the $\sup$-metric,
and  $B(C)\subset H^0_b(V_0, \calo_V)$ the set of all
functions $f$ with $|f|\leq C$. This is a normal family
(\ref{normal_family}), hence it is precompact in $C^0$-topology
(that is, in the topology of uniform convergence on compacts).
However, it is not compact in the $\sup$-topology, because
the closed ball in a Banach space is never compact 
(this statement is called ``the Riesz theorem'', see
\cite[Chapter 1, Theorem 4]{_Diestel_}.\footnote{Compare with the proof of 
	\cite[Theorem 2.14]{ov_indam}.} 
%or Exercise \ref{_Riesz_Exercise_}). 

Consider $(\gamma^k)^*:\; H^0_b(V_0, \calo_V)\arrow H^0_b(V_0, \calo_V)$
as a linear map of Banach spaces.
The map $(\gamma^k)^*$ takes any sequence $\{f_i\}$ converging in
$C^0$-topology to a sequence converging in the $\sup$-topology.
Indeed, $\{f_i\}$ uniformly converges on the closure
$\overline{\gamma^k(V_0)}$ which is compact because $\gamma^k(V_0) \Subset V_0$.
Since $(\gamma^k)^*(f_i)(z)= f_i(\gamma^k(z))$, the
function $f_i$ takes the same values on $\gamma^k(V_0)$ as
$(\gamma^k)^*(f_i)$ takes on $V_0$.
Therefore, the sequence $\{(\gamma^k)^*(f_i)\}$ uniformly converges
 on  $V_0$ for any $k>0$.

Since an open ball in $H^0_b(V_0, \calo_V)$ is precompact
in the $C^0$-topology, the operator
$(\gamma^k)^*:\; H^0_b(V_0, \calo_V)\arrow H^0_b(V_0, \calo_V)$
takes an open ball to a precompact set in the Banach space
$H^0_b(V_0, \calo_V)$. Therefore, the operator 
$(\gamma^k)^*$ is compact as a map of Banach spaces.

\hfill

{\bf Step 2:} 
A function $f \in H^0_b(V_0, \calo_V)$
is called {\bf $\gamma^k$-finite} if 
the sequence 
$\{f, (\gamma^k)^*(f), (\gamma^{2k})^*(f), 
(\gamma^{3k})^*(f),...\}$ belongs to a finite-dimensional
subspace of $H^0_b(V_0, \calo_V)$ (in the language used in
\ref{_F_finite_Definition_},
we would have to say {\bf $(\gamma^k)^*$-finite}).
Let $R$ be a compact injective operator
on a Banach space $B$. From \ref{_Riesz_Schauder_main_Theorem_}
(see also \ref{_root_space_RS_Remark_}), it follows that 
$B$ admits a Schauder basis%
\footnote{A {\bf Schauder basis} in a Banach space $W$
	is a set $\{x_i\}$ of vectors such that the closure of the
	space generated by $x_i$ is $W$, and the closure of the space generated 
	by all of $\{x_i\}$ except one of them does not contain the last one.}
such that the matrix for $R$ in this basis\index{function!$\gamma^*$-finite}
is a sum of finite-dimensional Jordan blocks. This implies
that the set of $\gamma^k$-finite vectors in $H^0_b(V_0, \calo_V)$
is dense (see also \ref{_finite_RS_Corollary_}).

\hfill

{\bf Step 3:} 
Every $\gamma^k$-finite function $f$ can be extended
to a holomorphic function on $V$. Indeed,
let $W_f$ be the finite-dimensional subspace
of $H^0_b(V_0, \calo_V)$ generated by 
$\{f, (\gamma^k)^*(f), (\gamma^{2k})^*(f), 
(\gamma^{3k})^*(f),...\}$. Since 
$(\gamma^k)^*$ takes $W_f$ to itself, and $W_f$ is finite-dimensional,
the map $(\gamma^k)^*$ is invertible on $W_f$.
However, $((\gamma^k)^*)^{-1}$ takes functions
defined on $V_0$ to functions defined on $\gamma^{-k}(V_0)$,
hence any $w\in W_f$ can be holomorphically extended to
$\gamma^{-k}(V_0)$. Repeating this argument, we extend
$w\in W_f$ from $V_0$ to the sets
$\gamma^{-mk}(V_0)\supset \gamma^{(-m+1)k}(V_0)\supset ... \supset 
\gamma^{-2k}(V_0)\supset \gamma^{-k}(V_0)\supset V_0$.
Since $\gamma$ is a contraction, we have
$\bigcup_m \gamma^{-mk}(V_0)=V$, hence all elements of $W_f$  
are holomorphically extended to $V$; this includes $f\in W_f$.

\hfill

{\bf Proof. Step 4:} 
It is clear that every  $\gamma$-finite function $f_1\in H^0(V, \calo_V)$ 
is $\gamma^k$-finite. Now we prove that, conversely,
every $\gamma^k$-finite function $f\in H^0(V, \calo_V)$ 
is in fact $\gamma$-finite. Indeed, let 
$W_f\subset H^0(V, \calo_V)$ be 
the subspace generated by 
$\{f, (\gamma^k)^*(f), (\gamma^{2k})^*(f), 
(\gamma^{3k})^*(f),...\}$, and $W'_f\supset W_f$ the subspace
generated by $\{f, (\gamma)^*(f), (\gamma^{2})^*(f), 
(\gamma^{3})^*(f),...\}$. 

Until the rest of this step, we
use $\sum_{i=0}^\infty z_i$ as a shorthand 
for the subspace generated by the
vectors $\{z_i\}$.
Since $\Z^{\geq 0}= \bigcup_{i=0}^{k-1} \big[i+(k\Z^{\geq 0})\big],$
the space
$W'_f= \sum_{i=0}^{k-1} (\gamma^i)^* \sum_{l\geq 0}^\infty  (\gamma^{lk})^*(f)$ 
is equal to the subspace generated by the finite-di\-men\-si\-o\-nal
spaces $W_f, \gamma^*(W_f), ..., (\gamma^{k-1})^*(W_f)$,
hence it is also finite-dimensional.

\hfill

{\bf Proof. Step 5:} 
Since the space $V$ is Stein, it admits
a closed holomorphic embedding to $\C^d$.
In other words, there exists a 
$d$-dimensional subspace $A\subset H^0(V, \calo_V)$ 
such that the corresponding map $V \arrow \C^d$
is a closed embedding.
Since the space of $(\gamma^k)$-finite functions is dense in 
$H^0_b(V_0, \calo_V)$, the space of $\gamma$-finite functions
is dense in $H^0(V, \calo_V)$ with $C^0$-topology. Let $P$ be a
precompact fundamental domain for the action of $\langle \gamma\rangle $ on $V$.
Since $H^0(V, \calo_V)^\gamma$ is $C^0$-dense in $H^0(V, \calo_V)$,
there exists a finite-dimensional
space $A_1= \C^d$ of $\gamma$-finite functions such that the
corresponding map $P \arrow \C^d$ is an embedding.
Let $W\subset H^0(V, \calo_V)$ be the smallest
$\gamma^*$-invariant subspace generated by $A_1$;
since all vectors in $A_1$ are $\gamma$-finite,
the space $W$ is finite-dimensional.
This makes the following diagram commutative
\begin{equation*}
	\begin{CD}
		V@>\Psi>> W \\
		@V{\gamma}VV  @VV{\gamma^*} V \\
		V @>\Psi >>  W
	\end{CD}
\end{equation*}
were $\Psi:\; V \arrow W$ is a closed embedding
of $V$ to $W=\C^N$. Removing the zero and taking
the quotient over $\gamma$, respectively $\gamma^*$, this gives a
closed holomorphic embedding 
$\frac{V\backslash x}{\langle \gamma\rangle}
\hookrightarrow \frac{W\backslash 0}{\langle \gamma^*\rangle}$,
where $\frac{W\backslash 0}{\langle \gamma^*\rangle}$ is by construction
a linear Hopf manifold.
\endproof

\hfill

\corollary\label{_LCK_pot_on_quatient_Stein_Corollary_}
Let $V$ be a Stein variety, and $\gamma:\; V \arrow V$ 
an invertible holomorphic contraction centered in $x$.
Assume that the group $\langle\gamma\rangle$ acts on $V\backslash x$ properly
discontinuously, and $V\backslash x$ is smooth. Then the quotient manifold
$\frac{V\backslash x}{\langle \gamma\rangle}$ admits
an LCK metric with potential.

\hfill

\proof 
By \ref{_Stein_by_contract_to_linear_Hopf_Theorem_},
$\frac{V\backslash x}{\langle \gamma\rangle}$ admits
a holomorphic embedding to a linear Hopf manifold $H$.
By \cite[Section 2.5]{ov_pams},
$H$ admits an LCK metric with potential.
Any closed complex submanifold of an LCK manifold with potential
is again LCK with potential, hence 
$\frac{V\backslash x}{\langle \gamma\rangle}$
is LCK with potential.
\endproof

%%%%%%%%%%%%%%%%%%%%%%%%%%%%%%%%%%%%%%%%%%%%%%%%%%%%%%%%%%%%
\section{Non-linear Hopf manifolds are LCK}
%%%%%%%%%%%%%%%%%%%%%%%%%%%%%%%%%%%%%%%%%%%%%%%%%%%%%%%%%%%%

%%%%%%%%%%%%%%%%%%%%%%%%%%%%%%%%%%%%%%%%%%%%%%%%%%%%%%%%%%%%
\definition\label{_Hopf_gene_Definition_}
Let $\gamma:\; \C^n \arrow \C^n$ be an invertible holomorphic
contraction centered in 0. Then the quotient
$\frac{\C^n \backslash 0}{\langle \gamma\rangle}$
is called {\bf a Hopf manifold}.\index{manifold!Hopf}

\hfill

\remark 
A Hopf manifold is diffeomorphic to $S^1 \times S^{2n-1}$.
In particular, it is never K\"ahler.

\hfill

The following corollary is directly implied
by \ref{_LCK_pot_on_quatient_Stein_Corollary_}.

\hfill

\corollary\label{_hopf_are_lck_corollary_}
Let $H$ be a Hopf manifold. Then $H$ admits an LCK
metric with potential.
\endproof

\hfill

\remark
Most of the results about the LCK manifolds with potential
rely upon the Andreotti-Siu and Rossi theorem used to construct
the Stein completion of the K\"ahler $\Z$-cover. This is why we
restrict most of the statements to manifolds of dimension $\geq 3$.
However, a Hopf manifold $H$ is covered by $\C^n \backslash 0$,
and this manifold already has the Stein completion $\C^n$.
Therefore, in the results about the Hopf manifolds
we don't need the restriction $\dim_\C H \geq 3$.

\hfill

\noindent{\bf Acknowledgment:} We are grateful to Florin Belgun for insightful discussions about the subject of this paper.

\hfill

\hfill

{\small
	
	\noindent {\sc Liviu Ornea\\
		University of Bucharest, Faculty of Mathematics and Informatics, \\14
		Academiei str., 70109 Bucharest, Romania}, and:\\
	{\sc Institute of Mathematics ``Simion Stoilow" of the Romanian
		Academy,\\
		21, Calea Grivitei Str.
		010702-Bucharest, Romania\\
		\tt lornea@fmi.unibuc.ro,   liviu.ornea@imar.ro}
	
	\hfill

	\noindent {\sc Misha Verbitsky\\
		{\sc Instituto Nacional de Matem\'atica Pura e
			Aplicada (IMPA) \\ Estrada Dona Castorina, 110\\
			Jardim Bot\^anico, CEP 22460-320\\
			Rio de Janeiro, RJ - Brasil }\\
		also:\\
		Laboratory of Algebraic Geometry, \\
		Faculty of Mathematics, National Research University 
		Higher School of Economics,
		6 Usacheva Str. Moscow, Russia}\\
	\tt verbit@verbit.ru, verbit@impa.br }

\end{document}